\documentclass[journal]{IEEEtran}

\usepackage{cite}

\usepackage{algorithmic}
\usepackage{algorithm2e}
\usepackage{caption}

\usepackage{color}
\usepackage{placeins}
\usepackage{float}
\usepackage{tabularx,colortbl}

\usepackage{subfigure}

\usepackage[cmex10]{amsmath}
\usepackage[pdftex]{graphicx}
\usepackage{amsfonts}

\usepackage{xcolor}

\usepackage{array}

\newcolumntype{L}[1]{>{\raggedright\let\newline\\\arraybackslash\hspace{0pt}}m{#1}}
\newcolumntype{C}[1]{>{\centering\let\newline\\\arraybackslash\hspace{0pt}}m{#1}}
\newcolumntype{R}[1]{>{\raggedleft\let\newline\\\arraybackslash\hspace{0pt}}m{#1}}

\newcommand{\mbold}[1]{\mathbf{#1}}

\newcommand{\mvar}[1]{\boldsymbol{#1}}

\begin{document}
%\title{A Capacity-based Matrix Design Method for Block Compressive Sensing Applications}
\title{Sensing Matrix Design via Capacity Maximization for Block Compressive Sensing Applications}

% use a single column layout for the different affiliations.
% In case the affiliation is too long then use \\ to create multiple lines, as shown in affiliation 4
% Too long lines can result in margin upload warnings ....
\author{\IEEEauthorblockN{Richard Obermeier and Jose~Angel~Martinez-Lorenzo}}% <-this % stops a space



\maketitle

\begin{abstract}
It is well established in the compressive sensing (CS) literature that sensing matrices whose elements are drawn from independent random distributions exhibit enhanced reconstruction capabilities. In many CS applications, such as electromagnetic imaging, practical limitations on the measurement system prevent one from generating sensing matrices in this fashion. Although one can usually randomized the measurements to some degree, these sensing matrices do not achieve the same reconstruction performance as the truly randomized sensing matrices. In this paper, we present a novel method, based upon capacity maximization, for designing sensing matrices with enhanced block-sparse signal reconstruction capabilities. Through several numerical examples, we demonstrate how our method significantly enhances reconstruction performance.
\end{abstract}

%\noindent \begin{abstract}
%Compressive Sensing (CS) theory states that sparse signals can be recovered from a small number of linear measurements $y=Ax$ using $\ell_1-$norm minimization techniques, provided that the sensing matrix satisfies a Restricted Isometry Property (RIP). Unfortunately, the RIP is difficult to verify in electromagnetic imaging applications, where the sensing matrix is computed deterministically. Although it provides weaker reconstruction guarantees than the RIP, the mutual coherence is a more practical metric for assessing the CS recovery properties of deterministic matrices. In this paper, we describe a method for minimizing the mutual coherence of sensing matrices in electromagnetic imaging applications. Numerical results for the design method are presented for a simple multiple monostatic imaging application, in which the sensor positions for each measurement serve as the design variables. These results demonstrate the algorithm's ability to both decrease the coherence and to generate sensing matrices with improved CS recovery capabilities.
%\end{abstract}


%{\smallskip \keywords antenna, propagation, measurement.}
% Keywords

% \fi
\IEEEpeerreviewmaketitle

%%%%%%%%%%%%%%%%%%%%%%%%%%%%%%%%%%%%%%%%%%%%%%%%%%%%%%
\vspace{7pt}
\section{Introduction}
\label{sec:intro}

% Undeter-determined linear systems ==> infinite number of solutions ==> need to regularize
A classical problem in science and engineering is reconstructing an unknown vector $\mbold{x} \in \mathbb{C}^N$ from a set of linear measurements $\mbold{y} = \mbold{A}\mbold{x} \in \mathbb{C}^M$. When $M < N$, there exist an infinite number of solutions satisfying $\mbold{y} = \mbold{A}\mbold{x}$ and so regularization techniques need to be employed in order to induce a unique solution. In practice, the regularization term is selected from prior knowledge of the unknown vector. When the vector is known to be sparse, then Compressive Sensing (CS) theory \cite{Candes2006,Candes2006a,Donoho2006} states that it can be recovered exactly as the solution to a convex and computationally tractable $\ell_1-$norm minimization problem, provided that the sensing matrix is ``well-behaved'' according to a performance metric such as the mutual coherence \cite{donoho2001} or the Restricted Isometry Property (RIP) \cite{candes2008}.

% CS ==> recover sparse vectors, usually done by applying norm-one approach, RIP and coherence

% Block CS ==> recover block-sparse vectors, extensions to RIP and coherence
CS theory also considers the case where the unknown vector is block sparse. When a signal is block sparse, the non-zero values are distributed over $K = N/L$ disjoint blocks of size $L$. Although block sparse signals can be reconstructed using the standard techniques, such as $\ell_1-$norm minimization, that are applied to general sparse signals, specialized techniques based on joint $\ell_2 / \ell_1$ minimization have been shown to provide better reconstruction performance  \cite{stojnic2010ell_,zeinalkhani2015iterative,garg2011block,eldar2010average,eldar2009robust,eldar2010block}. Unsurprisingly, extensions to the coherence \cite{eldar2010block} and RIP \cite{stojnic2009reconstruction,gao2017new,eldar2009robust,eldar2010average} determine whether or not the block sparse recovery techniques will be successful for a given sensing matrix.

% Can be difficult to design good sensing matrices that satisfy RIP / coherence criteria. List some existing methods. Briefly describe their shortcomings and introduce the capacity-based method
In general, one cannot deterministically generate sensing matrices that satisfy the RIP or the block-sparse variant. Often times, researchers will resort to random matrix theory in order to generate sensing matrices that satisfy the RIP with high probability. Unfortunately, this approach cannot be used in many applications, such as electromagnetic imaging, where the elements of the sensing matrix are constrained by practical limitations. In this paper, we introduce a method based upon maximizing the sensing capacity \cite{migliore2008electromagnetics,Migliore2011} for designing sensing matrices with enhanced block sparse signal reconstruction capabilities.

% Paper organization
The remainder of this paper is organized as follows. In Section \ref{sec:motivation}, we discuss the motivation for using the sensing capacity as the design metric. In Section \ref{sec:design_method}, we formulate the capacity-based design method and describe how it can be solved using the method of multipliers \cite{nocedal2006numerical}. In Section \ref{sec:prev_work}, we discuss previous work that has been performed on designing sensing matrices for block sparse signal reconstruction problems. In Section \ref{sec:results}, we present results for several design scenarios to demonstrate the effectiveness of the algorithm. Finally, in Section \ref{sec:conclusion} we conclude the paper by describing several other applications where the design algorithm can be applied.

%%%%%%%%%%%%%%%%%%%%%%%%%%%%%%%%%%%%%%%%%%%%%%%%%%%%%%
\vspace{7pt}
\section{Motivation}
\label{sec:motivation}
Consider the noise-corrupted linear system $\mbold{y} = \mbold{A}\mbold{x} + \mbold{n}$, where $\mbold{x} \in \mathbb{C}^N$, $\mbold{y},\mbold{n} \in \mathbb{C}^M$, $\mbold{A} \in \mathbb{C}^{M\times N}$, and $M < N$. It is assumed here that $\mbold{A}$ has normalized columns. Suppose that the unknown vector is known to be block sparse with block size $L$ and let us denote $\mbold{P}_k \in \{0,1\}^{L \times N}, ~k = 1, \ldots, K$ as the binary projection matrix that extracts the elements of $\mbold{x}$ in the $k-$th block. Note that because $\mbold{P}_k$ is a projection matrix, $\mbold{P}_k\mbold{P}_k^T = \mbold{I}_{L,L}$, the identity matrix, and $\mbold{P}_k\mbold{P}_j^T = \mbold{0}_{L,L}$, the zero matrix, for $k\neq j$. In order to induce sparsity in the solution vector, one would ideally use a mixed $\ell_2 / \ell_0$ objective function, where the $\ell_0$-``norm'' simply counts the number of non-zeros. Unfortunately, this problem is $NP$-hard, and so it cannot be solved even for moderately sized problems. However, block CS theory states that the vector can be stably recovered using the following joint $\ell_2 / \ell_1$ technique  \cite{stojnic2010ell_,zeinalkhani2015iterative,garg2011block,eldar2010average,eldar2009robust,eldar2010block}
\begin{alignat}{3}
& \underset{\mvar{x}}{\text{minimize }}~~ & &\sum_{k=1}^K \|\mbold{P}_k\mvar{x}\|_{\ell_2}  \label{eq:block_cs} \\
&\text{subject to } &  &\|\mbold{A}\mvar{x} - \mbold{y}\|_{\ell_2} \le \eta \nonumber
\end{alignat}
provided that the sensing matrix $\mbold{A}$ is ``well-behaved'' according to some design metric. The most powerful design metric is the block RIP \cite{stojnic2009reconstruction,gao2017new,eldar2009robust,eldar2010average}, which can be defined as follows. For a fixed block sparsity level $T$, the block restricted isometry constant $\delta_{L,T}$ is the smallest positive constant such that
 \begin{align}
 (1-\delta_{L,T})\|\mbold{x}_i\|_{\ell_2}^2 \le \|\mbold{A}_i\mbold{x}_i\|_{\ell_2}^2 \le (1+\delta_{L,T})\|\mbold{x}_i\|_{\ell_2}^2 \label{eq:block_ric}
 \end{align}
where $\mbold{x}_i = \mbold{\Phi}_i \mbold{x}$ and $\mbold{A}_i = \mbold{A}\mbold{\Phi}_i^T$, for all projection matrices $\mbold{\Phi}_i \in \{0,1\}^{LT \times N},~i = 1,\ldots,\binom{K}{T}$ obtained by concatenating $T$ of the $K$ projection matrices. Generally speaking, the block RIP requires $\delta_{L,T}$ to be small. Note that if $L=1$, then Eq. \ref{eq:block_cs} reduces to $\ell_1$-norm minimization and Eq. \ref{eq:block_ric} reduces to the standard RIP. 

% Formally describe the RIP and how it extends to block CS applications. Cite papers.

% Tie together the relationship between RIP and sensing capacity. Cite papers.
The block RIP can also be analyzed from the perspective of information theory. The $\epsilon$-capacity \cite{migliore2008electromagnetics,Migliore2011}, also referred to as the sensing capacity, determines the amount of information that can be transmitted through a linear mapping within an uncertainty level $\epsilon$. For the linear mapping $\mbold{A}_i$, the sensing capacity can be expressed as follows:
\begin{align}
H_\epsilon(\mbold{A}_i) = \frac{1}{2}\log_2\left(\det \mbold{A}_i^H\mbold{A}_i\right) = \sum_{m=1}^{LT}\log_2\left(\frac{\sigma_{m,i}}{\epsilon}\right) \label{eq:cap}
\end{align}
where $\sigma_{m,t}$ is the $m-$th singular value of $\mbold{A}_i$. From the definition of the block RIP, it is easy to show that the sensing capacity $H_\epsilon(\mbold{A}_i)$ is bounded by:
\begin{align}
H_\epsilon(\mbold{A}_i)  \ge \frac{LT}{2} \left(\log_2\left(\frac{\sqrt{1-\delta_{L,T}}}{\epsilon}\right)+\log_2\left(\frac{\sqrt{1+\delta_{L,T}}}{\epsilon}\right)\right) \label{eq:cap_bound}
\end{align}
which implies that $\delta_{L,T}$ determines the minimum amount of information that can be transmitted by $T$ block block sparse vectors using the linear mapping $\mbold{y} = \mbold{A}\mbold{x}$. 

% Propose the capacity-based approach for a subset of groups of columns. Feasible for block CS applications.
The motivation for using the sensing capacity as the design metric is apparent at this point. Maximizing the minimum capacity over a set of sub-matrices $\mbold{A}_i$ allows us to decrease the block restricted isometry constant $\delta_{L,T}$. Although it is impractical to optimize over all $\binom{K}{T}$ sub-matrices for our desired block sparsity level $T$, we shall see in Section \ref{sec:results} that we can obtain significant improvements by optimizing over $\binom{K}{2} = \frac{K(K-1)}{2}$ sub-matrices of size $M \times 2L$. This approximation is analogous to that taken by mutual coherence minimization, which minimizes the restricted isometry constant $\delta_2$. In fact, when $L=1$ and $T=2$, Eq. \ref{eq:cap_bound} is a tight bound.

%%%%%%%%%%%%%%%%%%%%%%%%%%%%%%%%%%%%%%%%%%%%%%%%%%%%%%
\vspace{7pt}
\section{Capacity-based Design Method}
\label{sec:design_method}

% Describe the method and some of the implementation details. Perhaps refer the reader to the coherence minimization paper for some of the details.
Suppose that the sensing matrix $\mbold{A} \in \mathbb{C}^{M \times N}$ is a function of $\mbold{p} \in \mathbb{C}^{P}$ design variables according to the nonlinear and differentiable relationship $\mbold{A} = \mbold{F}(\mbold{p})$. Without loss of generality, we will assume that this function outputs the sensing matrix with normalized columns. If necessary, one can replace the function $\mbold{f}_m(\mbold{p})$, which computes the $m-$th column of the sensing matrix, with $\hat{\mbold{f}}_m(\mbold{p}) = \frac{\mbold{f}_m(\mbold{p})}{\|\mbold{f}_m(\mbold{p})\|_{\ell_2}}$, which is differentiable everywhere except $\mbold{f}_m(\mbold{p}) = \mbold{0}_M$. Following the notation from the previous section, we define the projection matrices $\mbold{\Phi}_r \in \{0,1\}^{N \times M_r},~r=1,\ldots,R$ for the $R$ blocks on which the capacity will be evaluated. The design algorithm then seeks the minimizer to the following non-convex optimization program:
\begin{alignat}{3}
& \underset{\mbold{p}}{\text{minimize }}~~ & & \underset{r=1,\ldots,R}{\text{max }} -\log\det\left(\mbold{\Phi}^T_{r}\mbold{F}^H(\mbold{p})\mbold{F}(\mbold{p})\mbold{\Phi}_{r}  + \beta \mathbf{I}_{M_r,M_r}\right) \nonumber \\
&\text{subject to } &  & \mathbf{p}\in Q_p \label{eq:design_alg}
\end{alignat}
where $\beta$ is a small positive constant that ensures that the arguments to $\det$ are positive-definite, and $Q_p$ is the feasible set for the design variables. In other words, this optimization program seeks the design variables $\mathbf{p}$ that maximizes the smallest capacity of the sub-matrices $\mbold{F}(\mbold{p})\mbold{\Phi}_{r},~r=1,\ldots,R$.

Eq. \ref{eq:design_alg} can be solved using the method of multipliers \cite{nocedal2006numerical} with some modifications. To start, we introduce the auxiliary variable $\mbold{c} = \left(c_1, \ldots, c_R\right)^T \in \mathbb{R}^R$ to represent the capacities of the sub-matrices. With this modification, the capacity optimization problem can be expressed in the following equivalent form:
\begin{alignat}{3}
& \underset{\mbold{p}, c_1,\ldots,c_R}{\text{minimize }}~~ & & \underset{r=1,\ldots,R}{\text{max }} -c_r \label{eq:design_alg2} \\
&\text{subject to } &  & \mathbf{p}\in Q_p \nonumber \\
& & & c_r = \log\det\left(\mbold{\Phi}^T_{r}\mbold{F}^H(\mbold{p})\mbold{F}(\mbold{p})\mbold{\Phi}_{r} + \beta \mathbf{I}_{M_r,M_r}\right)  \nonumber
\end{alignat}
We can make one more modification to this problem due to the following observation: because the sub-matrices $\mbold{F}(\mbold{p})\mbold{\Phi}_{r}$ have normalized columns, $\operatorname{tr}(\mbold{\Phi}^T_{r}\mbold{F}^H(\mbold{p})\mbold{F}(\mbold{p})\mbold{\Phi}_{r}) = M_r$ and $\log\det\left(\mbold{\Phi}^T_{r}\mbold{F}^H(\mbold{p})\mbold{F}(\mbold{p})\mbold{\Phi}_{r}\right) \le 0$, where equality holds only when $\mbold{\Phi}^T_{r}\mbold{F}^H\mbold{F}(\mbold{p})\mbold{\Phi}_{r} = \mathbf{I}_{M_r,M_r}$. As a result, we can replace $\underset{r=1,\ldots,R}{\text{max }} -c_r$ in Eq. \ref{eq:design_alg2} with $\|\mbold{c}\|_{\ell_\infty}$. This allows us to utilize the proximal operator for the $\ell_\infty$-norm in the optimization procedure, so that we instead solve the following optimization program: 
\begin{alignat}{3}
& \underset{\mbold{p}, \mbold{c}}{\text{minimize }}~~ & & \|\mbold{c}\|_{\ell_\infty} \label{eq:design_alg3} \\
&\text{subject to } &  & \mathbf{p}\in Q_p \nonumber \\
& & & c_r = \log\det\left(\mbold{\Phi}^T_{r}\mbold{F}^H(\mbold{p})\mbold{F}(\mbold{p})\mbold{\Phi}_{r} + \beta \mathbf{I}_{M_r,M_r}\right)  \nonumber
\end{alignat}
Eq. \ref{eq:design_alg2} and \ref{eq:design_alg3} are only equivalent when $\mbold{F}(p)$ has normalized columns. In actuality, Eq. \ref{eq:design_alg3} minimizes the maximum absolute value of the capacity. So, when $\mbold{F}(p)$ does not have normalized columns, Eq. \ref{eq:design_alg3} will drive the capacities closer to zero instead of maximizing the smallest capacity. 

Eq. \ref{eq:design_alg3} has a very similar form to the coherence minimization algorithm displayed in Eq. 11 of \cite{obermeier2017coherence}: simply replace the coherence equality constraints with the capacity equality constraints. As a result, the method of multipliers \cite{nocedal2006numerical} approach described in \cite{obermeier2017coherence} can also be used to solve Eq. \ref{eq:design_alg3}, provided that the feasible set $Q_p$ has an easy to compute proximal operator. Formally, the scaled Augmented Lagrangian can be written as follows:
\begin{align}
\mathcal{L}_{\mathcal{A}}(&\mbold{p},\mbold{c},\boldsymbol{\gamma}; \rho) = \|\mbold{c}\|_{\ell_\infty} + I_{Q_p}(\mbold{p}) +  \label{eq:auglag1} \\
&\sum_{r=1}^R \frac{\rho}{2}\left|c_{r} - \log\det\left(\mbold{\Phi}^T_{r}\mbold{F}^H(\mbold{p})\mbold{F}(\mbold{p})\mbold{\Phi}_{r} + \beta \mathbf{I}_{M_r,M_r}\right) + \gamma_r/\rho \right|^2 \nonumber
\end{align}
where $\boldsymbol{\gamma} \in \mathbb{R}^R$ are the Lagrange multipliers. The method of multipliers solves Eq. \ref{eq:design_alg3} by solving a series of unconstrained problems of the form of Eq. \ref{eq:auglag1}, where $\boldsymbol{\gamma}$ is held fixed. The unconstrained sub-problems can be solved using an alternating minimization schema, in which $\mbold{c}$ is updated by evaluating the proximal operator for the $\ell_\infty$-norm and $\mbold{p}$ is updated using a proximal gradient update. For details, the reader is referred to \cite{obermeier2017coherence}. When a given instance of Eq. \ref{eq:auglag1} is solved, the Lagrange multipliers are updated as follows:
\begin{align}
\gamma_{r}^{(k+1)} =~& \gamma_{r}^{(k)} + \rho^{(k)}\bigg(c_{r}^{(k)} - \\
&\log\det\left(\mbold{\Phi}^T_{r}\mbold{F}^H(\mbold{p})\mbold{F}(\mbold{p})\mbold{\Phi}_{r} + \beta \mathbf{I}_{M_r,M_r}\right)\bigg) \nonumber
\end{align}
where the superscripts indicate the iteration number, i.e. the Lagrange multiplier $\gamma_{r}^{(k)}$ is used on the $k-$th instance of Eq. \ref{eq:auglag1}. To ensure that the algorithm converges to a stationary point, it is often necessary to increase $\rho$ at each iteration. Our design method utilizes the update approach described in \cite{nocedal2006numerical}. The optimization procedure is summarized in Algorithm \ref{alg:auglagalg}. 
\begin{algorithm}
\caption{Summary of the Augmented Lagrangian update procedure for the capacity maximization problem of Eq. \ref{eq:design_alg3}}
\label{alg:auglagalg}
Choose the initial values for $\mbold{p}^{(0)}$, $\rho^{(1)}$ \;
Set $c_{r}^{(0)} = \log\det\left(\mbold{\Phi}^T_{r}\mbold{F}^H(\mbold{p}^{(0)})\mbold{F}(\mbold{p}^{(0)})\mbold{\Phi}_{r} + \beta \mathbf{I}_{M_r,M_r}\right)$, $\gamma_{r}^{(1)} = 0$ \;
\For{k = 1,2,3\ldots}{
 Solve the unconstrained subproblem
\[
 \left(\mbold{p}^{(k)},\mbold{c}^{(k)}\right) = \underset{\mbold{p},\mbold{u}}{\operatorname{argmin }} ~\mathcal{L}_{\mathcal{A}}\left(\mbold{p},\mbold{u},\boldsymbol{\gamma}^{(k)};\rho^{(k)}\right)
\]
\nl
  Update the dual variables
\begin{align}
\gamma_{r}^{(k+1)} =~& \gamma_{r}^{(k)} + \rho^{(k)}\bigg(c_{r}^{(k)} - \nonumber \\
&\log\det\left(\mbold{\Phi}^T_{r}\mbold{F}^H(\mbold{p})\mbold{F}(\mbold{p})\mbold{\Phi}_{r} + \beta \mathbf{I}_{M_r,M_r}\right)\bigg) \nonumber
\end{align}
\normalsize
\nl
 Compute $\rho^{(k+1)}$ using the method described in \cite{nocedal2006numerical}
}
\end{algorithm}

%To handle the un-normalized scenario, one can instead introduce an auxiliary variable $t$ and slack variables $\mbold{s} = \left(s_1, s_2, \ldots, s_R\right)^T$ to Eq. \ref{eq:design_alg2} and solve the following optimization program:
% \begin{alignat}{3}
%& \underset{\mbold{p}, \mbold{s}, t}{\text{minimize }}~~ & & t \label{eq:design_alg4} \\
%&\text{subject to } &  & \mathbf{p}\in Q_p \nonumber \\
%& & & t +  s_r + \log\det\left(\mbold{\Phi}^T_{r}\mbold{F}^H(\mbold{p})\mbold{F}(\mbold{p})\mbold{\Phi}_{r} + \beta \mathbf{I}_{M_r,M_r}\right) = 0 \nonumber \\
%& & & s_r \ge 0 \nonumber
%\end{alignat}

%\cite{obermeier2017coherence}

%%%%%%%%%%%%%%%%%%%%%%%%%%%%%%%%%%%%%%%%%%%%%%%%%%%%%%
\vspace{7pt}
\section{Comparison with Previous Work}
\label{sec:prev_work}

% Talk about the projection-based / subblock coherence methods. General shortcoming: can only be applied to sensing matrices that are linear functions of the design variables (projection matrix).

In practice, it is not possible to deterministically design a sensing matrix that satisfies the RIP or the block-sparse variant, and so researchers have instead focused on the coherence-based metric \cite{zelnik2011sensing,li2013projection,li2016block,qin2017novel}. These methods are limited in that they can only be applied to sensing matrices that are linear projections of a dictionary, i.e. $\mbold{A} = \mbold{\Phi} \mbold{D}$, where the elements of $\mbold{\Phi}$ are the design variables. Our method can be used to design sensing matrices that are nonlinear functions of the design variables, provided that the relationship is differentiable over the feasible set. In addition, our method can more directly optimize the block restricted isometry constants. As we mentioned in Section \ref{sec:motivation}, it is reasonable to optimize the capacity of the $R=\binom{K}{2} = \frac{K(K-1)}{2}$ sub-matrices of size $M \times 2L$, which indirectly optimizes the block restricted isometry constant $\delta_{L,2}$. 

% Capacity-based method takes a different approach and can be used for sensing matrices that are non-linear functions of the design variables.

%%%%%%%%%%%%%%%%%%%%%%%%%%%%%%%%%%%%%%%%%%%%%%%%%%%%%%%%%
\vspace{7pt}
\section{Numerical Results}
\label{sec:results}

% 1) Pulse reconstruction from incomplete fourier measurements
\subsection{Pulse Reconstruction Problem}
In the first example, the design algorithm was applied to a pulse reconstruction problem. Consider the scenario where we wish to reconstruct a time-series signal from a set of incomplete Fourier measurements. Formally, the $m-$th measurement can be expressed as follows:
\begin{align}
y_m = \sum_{n=1}^N x_n e^{-\jmath\omega_m n}
\end{align}
where $\omega_m$ is the digital frequency of the $m-$th measurement. The pulses were known to be distributed on $K=16$ non-overlapping segments of a fixed width $L=32$ samples. This is a simplified example of a sparsely used Time Division Multiple Access (TDMA) communication network. The objective, then, was to select the measurement frequencies $\omega_m$ such that the minimum capacities over all $\binom{16}{2} = 120$ pairs of blocks was maximized according to the design parameters and constraints displayed in Table \ref{tab:tab0}. Note that, due to the modulo nature of $\omega_m$, it was considered unbounded.
\begin{table}[h!]
\centering
\begin{tabular}{ |C{1.3cm}||C{3.1cm}||C{2.2cm}|  }
 \hline
 \multicolumn{3}{|c|}{Design Parameters and Constraints} \\
 \hline
 Parameter & Baseline Value & Constraint \\
 \hline
$M$ & $256$ & $-$ \\
$N$ & $512$ & $-$ \\
$K$ & $16$ & $-$ \\
$L$ & $32$ & $-$ \\
$\omega_m$ & Randomly distributed between $-\pi$ and $\pi$ & Unbounded \\
%$k$ & $62.8319 \text{m}^{-1}$ & $k = 62.8319 \text{m}^{-1}$ \\
 \hline
\end{tabular}
\caption{Summary of design parameters and constraints for the pulse reconstruction sensing matrix design problem.}
\label{tab:tab0}
\end{table}

\begin{figure}[h!]
\centering
\includegraphics[width=8.375cm, clip=true]{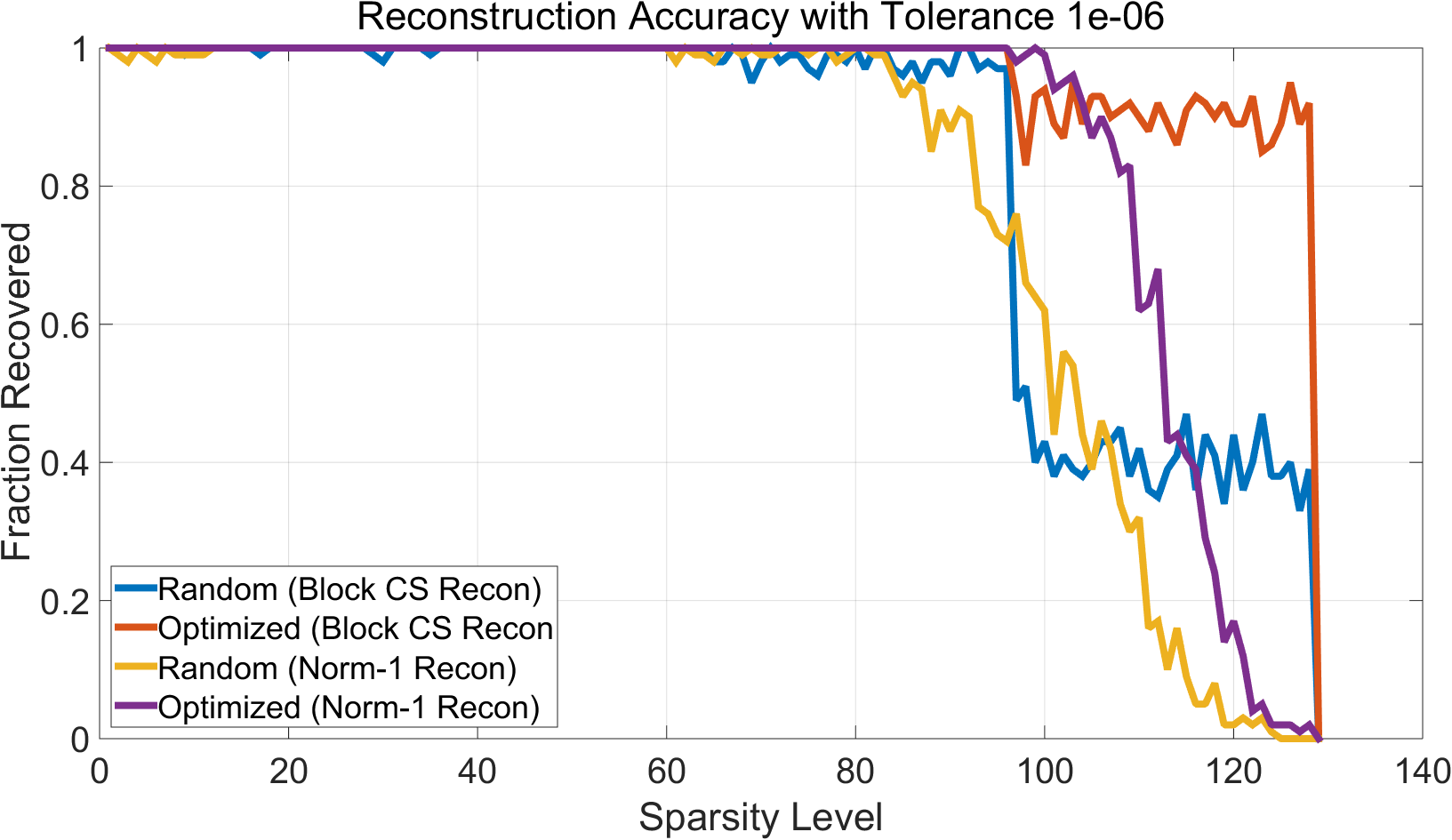}
    \caption{Numerical comparison of the reconstruction accuracies of joint $\ell_2/\ell_1$ reconstruction and standard $\ell_1$ reconstruction using the baseline random and optimized designs for the pulse reconstruction problem.}
   \label{fig:tdma_results}
\end{figure}

The capacity of the optimized design was increased from $-10.7$ to $-5.0$. While this may look like only a modest improvement at first glance, this had a significant affect on reconstruction performance. Figure \ref{fig:tdma_results} displays the CS reconstruction accuracies achieved by the baseline and optimized sensing matrices when joint $\ell_2/\ell_1$ and standard $\ell_1$ reconstruction techniques are used. These results were generated by reconstructing $100$ vectors at each sparsity level $S=1,\ldots,M/2$ (block sparsity $S_B = \frac{S}{L}$) and comparing the solutions to the ground truth vectors. The $\ell_1$-norm minimization results were included to provide a comparison with the joint $\ell_2/\ell_1$ results. Unsurprisingly, joint $\ell_2/\ell_1$ minimization outperformed $\ell_1$ minimization for each of the sensing matrices. Remarkably, the optimized sensing matrix was able to reconstruct $>90\%$ of block-sparse vectors up to a block sparsity $S_B=4$ (total sparsity $S=128$) using joint $\ell_2/\ell_1$ minimization, whereas the baseline random sensing matrix reconstructed $<50\%$. It is important to note that exact reconstruction cannot be guaranteed for total sparsity levels greater than $M/2$ ($128$ for this problem). Although it does not achieve the theoretical limit, the optimized sensing matrix achieves a level of performance that is significantly better than that of the randomized sensing matrix.

% 2) 2D Electromagnetic Imaging Problem
\subsection{Electromagnetic Imaging Problem}
\label{sec:emag_prob}
The design algorithm was also applied to an electromagnetic imaging application, in which a single transmitting and receiving antenna was used to excite a region of interest with a single frequency. The discretized measurement process for this system can be modeled as follows:
\begin{align}
y_m = \sum_{n=1}^N x_n e^{-j 2k \|\mbold{r}_m - \mbold{r}_n\|_{\ell_2}} = \sum_{n=1}^N A_{mn} x_n \label{eq:sensing_model}
\end{align}
where $y_m$ is the $m-$th scattered field measurement, $\mbold{r}_m$ is the position of the $m-$th antenna, $\mbold{r}_n$ is the $n-$th position in the imaging region, $k$ is the wavenumber, and $x_n$ is the reflectivity at the $n-$th position in the imaging region. Keeping the wavenumber fixed, the objective was to select the antenna positions $\mbold{r}_m$. 
\begin{table}[h!]
\centering
\begin{tabular}{ |C{1.3cm}||C{3.1cm}||C{2.2cm}|  }
 \hline
 \multicolumn{3}{|c|}{Design Parameters and Constraints} \\
 \hline
 Parameter & Baseline Value & Constraint \\
 \hline
$M$ & $64$ & $-$ \\
$N$ & $144$ & $-$ \\
$K$ & $9$ & $-$ \\
$L$ & $16$ & $-$ \\
$\mbold{r}_n$  & $5\lambda$ by $5\lambda$ grid centered at origin & $-$\\
$\mbold{r}_m$ & Uniformly spaced over $5\lambda$ by $5\lambda$ grid at $z=5\lambda$ & $|x_m| \le 2.5\lambda$  $|y_m| \le 2.5\lambda$ $z_m = 5\lambda$\\
%$k$ & $62.8319 \text{m}^{-1}$ & $k = 62.8319 \text{m}^{-1}$ \\
 \hline
\end{tabular}
\caption{Summary of design parameters and constraints for the electromagnetic imaging sensing matrix design problem.}
\label{tab:tab1}
\end{table}

\begin{figure}[h!]
\centering
\includegraphics[width=6.5cm, clip=true]{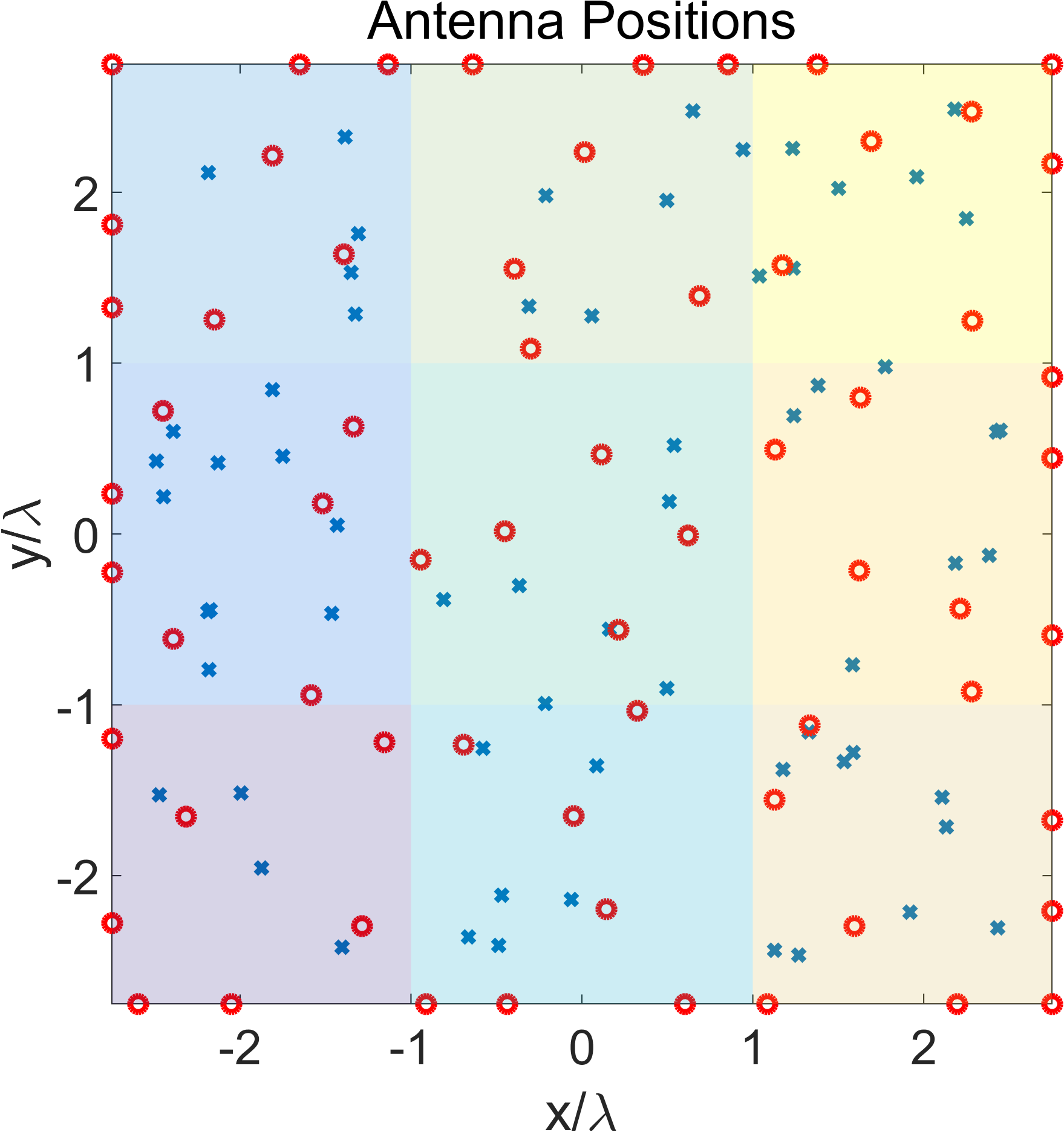}
    \caption{Antenna positions of the baseline (blue) and optimized (red) designs. The shaded boxes in the background represent the squares on which the capacity was evaluated.}
   \label{fig:ant_pos}
\end{figure}
\begin{figure}[h!]
\centering
\includegraphics[width=8.375cm, clip=true]{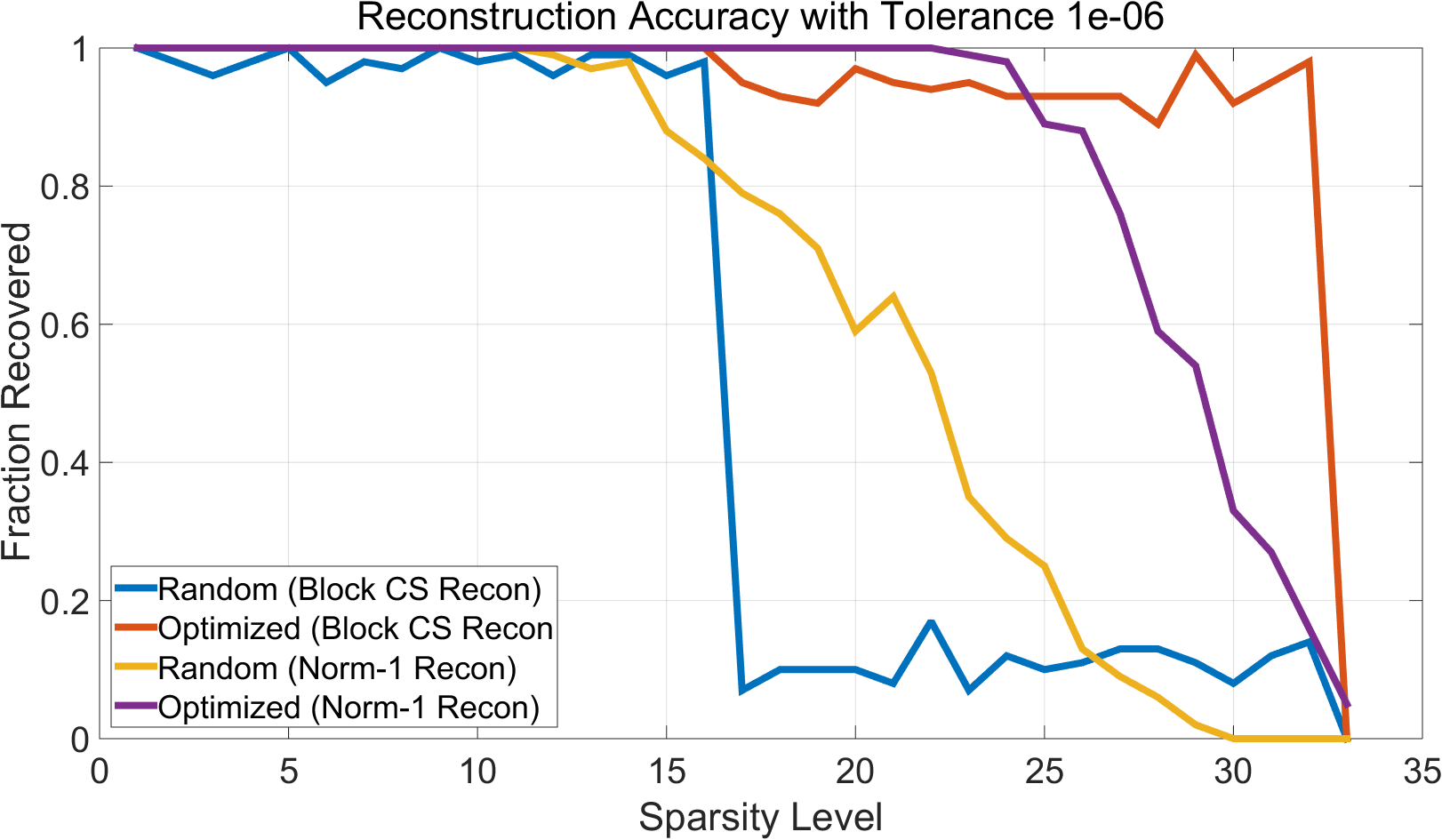}
    \caption{Numerical comparison of the reconstruction accuracies of joint $\ell_2/\ell_1$ reconstruction and standard $\ell_1$ reconstruction using the baseline random and optimized designs for the electromagnetic imaging problem.}
   \label{fig:results}
\end{figure}

\begin{figure}[h!]
\centering
\includegraphics[width=6.5cm, clip=true]{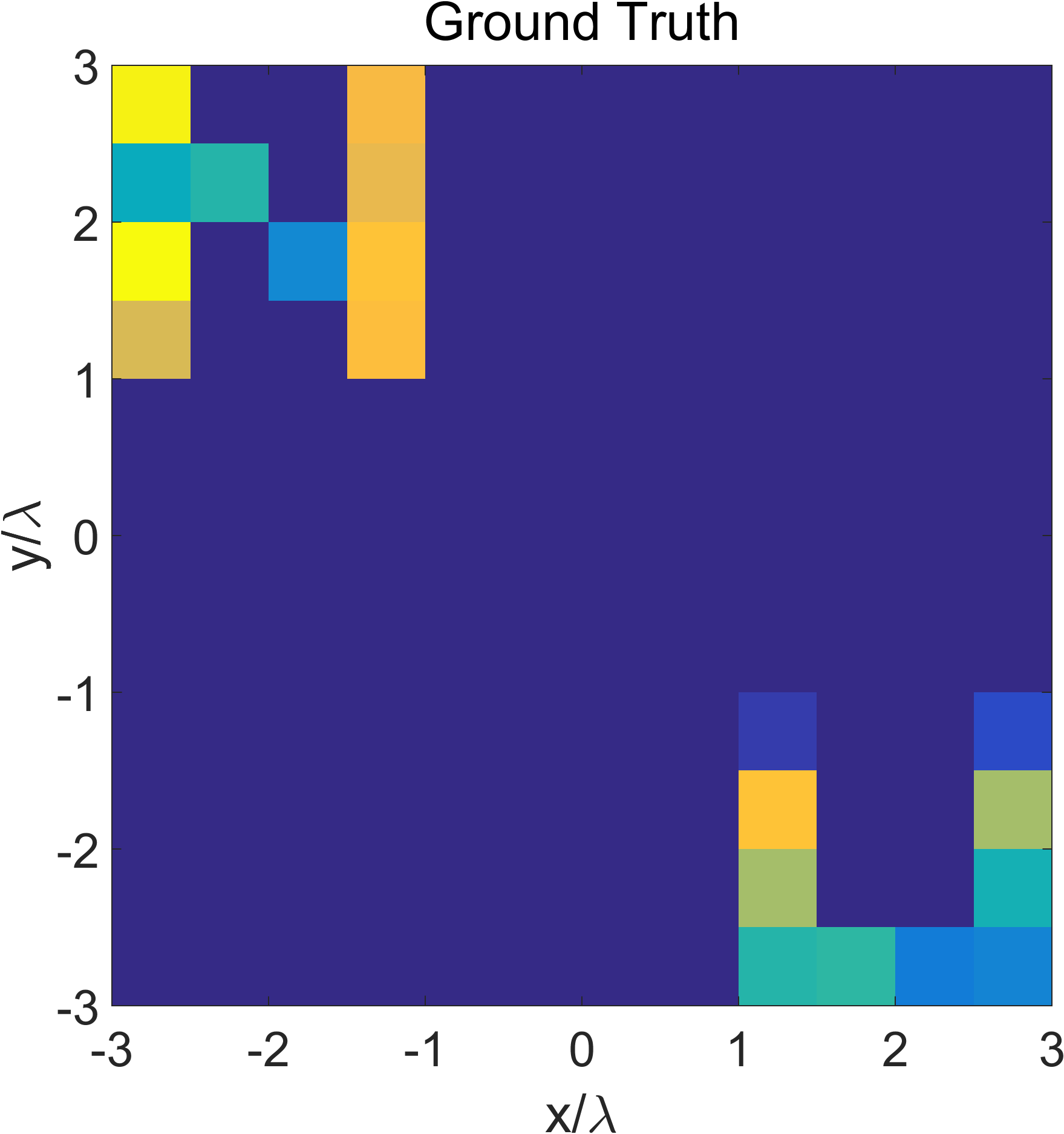}
    \caption{Magnitude of the ground-truth reflectivity}
   \label{fig:im_gt}
\end{figure}

\begin{figure}[h!]
\centering
\includegraphics[width=6.5cm, clip=true]{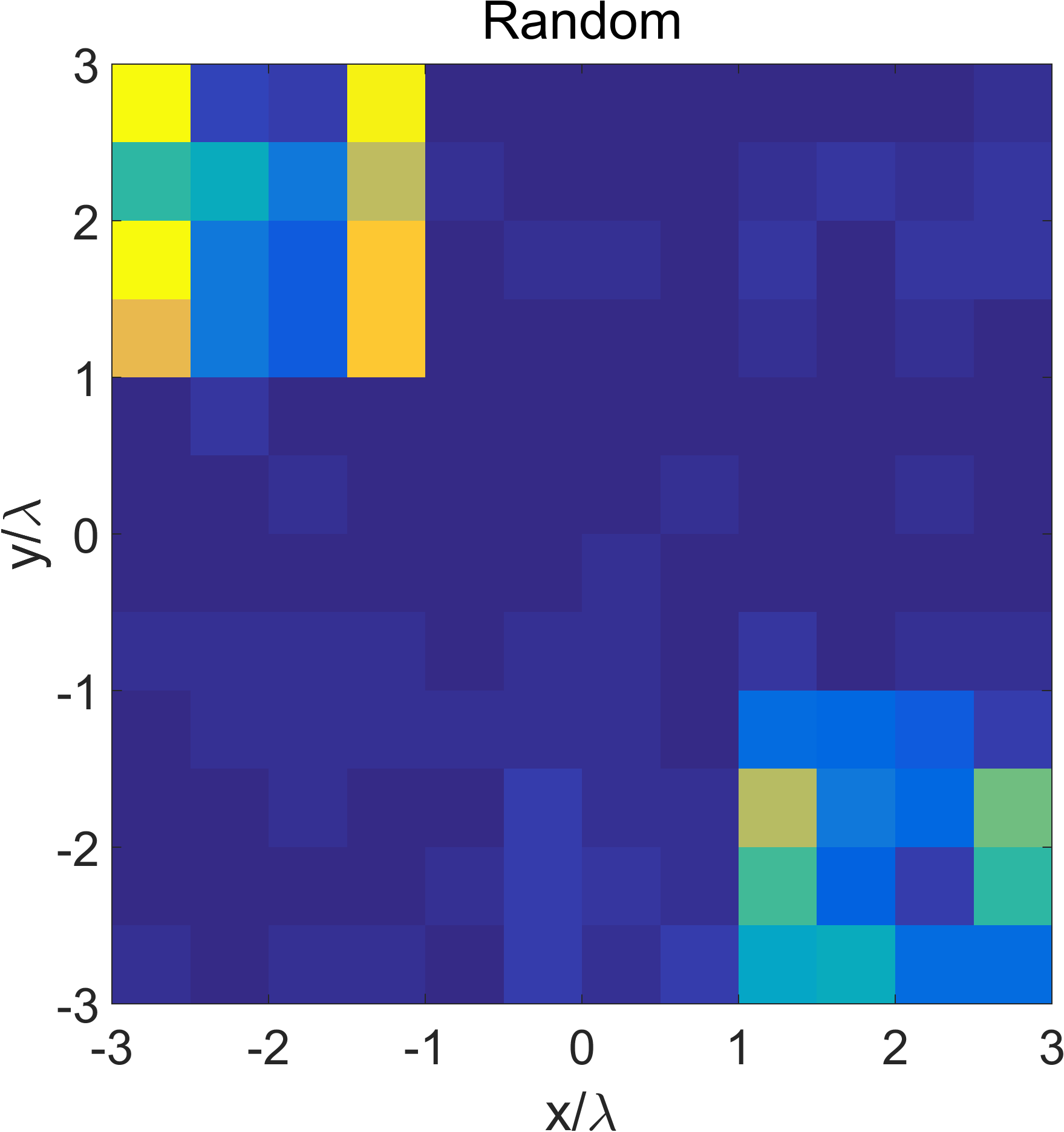}
    \caption{Magnitude of the reflectivity reconstructed by the baseline random sensing matrix using joint $\ell_2/\ell_1$ minimization. Normalized error = $0.2863$.}
   \label{fig:im_random}
\end{figure}

\begin{figure}[h!]
\centering
\includegraphics[width=6.5cm, clip=true]{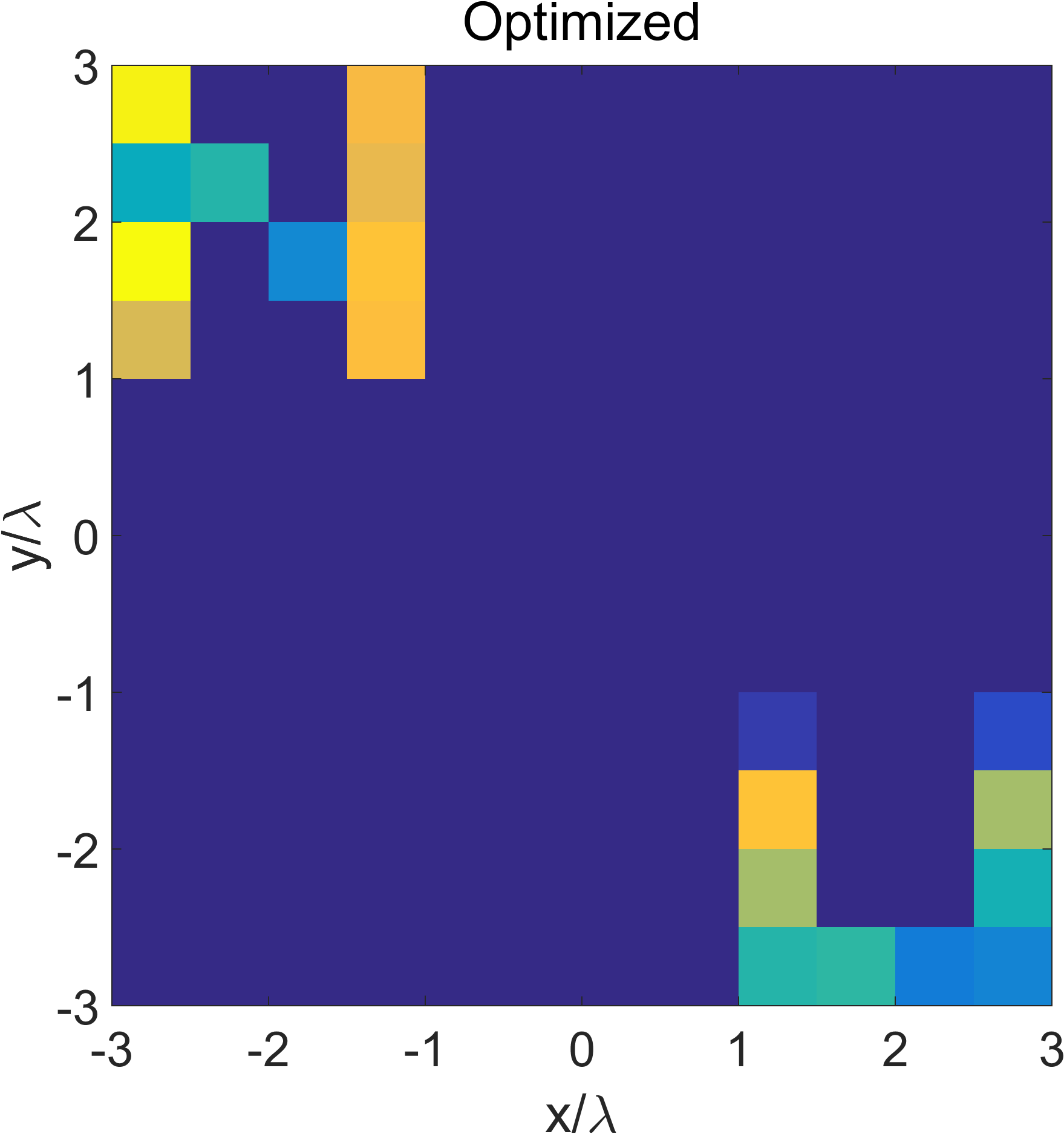}
    \caption{Magnitude of the reflectivity reconstructed by the baseline random sensing matrix using joint $\ell_2/\ell_1$ minimization. Normalized error = $0.0$.}
   \label{fig:im_opt}
\end{figure}

Table \ref{tab:tab1} displays the design parameters and constraints for the optimization problem, and Figure \ref{fig:ant_pos} displays the positions of the baseline random antenna configuration, which was used as the starting point to the optimization procedure, and the positions of the optimized antenna configuration. The shaded blocks in the background of Figure \ref{fig:ant_pos} represent the nine blocks on which the unknown signal was known to be block-sparse. The optimization procedure was configured so that the minimum capacity of all $36$ pairs of blocks was maximized. The optimized design achieved a minimum capacity of $-3.3$, which is a significant improvement over the baseline design, $-12.6$. This directly led to an improvement in CS reconstruction accuracy, as can be seen in Figure \ref{fig:results}. The optimized antenna positions were able to reconstruct $>90\%$ of block-sparse vectors up to a block sparsity $S_B=2$ (total sparsity $S=32$), whereas the baseline random positions reconstructed $<20\%$. Once again, the fact that the optimized design performs so well up to the theoretical maximum sparsity level, $M/2=32$, truly demonstrates the capabilities of the design method. A specific instance of the planar reconstruction problem is displayed in Figures \ref{fig:im_gt} - \ref{fig:im_opt}, which display the ground-truth reflectivity, the reflectivity reconstructed by the baseline random sensing matrix, and the reflectivity reconstructed by the optimized sensing matrix.


% 3) General sensing matrix problem. Benefits aren't as good in this example, which is somewhat expected.
\subsection{General Linear System}
In the final example, the design algorithm was tested against a general linear system. $\mbold{y} = \mbold{A}\mbold{x}$. The objective for this problem was to optimize the $MN$ coefficients $a_{mn} \in \mathbb{C}$ of the sensing matrix. This is the ideal design scenario, as we have complete control over the sensing matrix. The initial values for the sensing matrix were set by drawing values from i.i.d. complex Normal distributions. The full set of design parameters and constraints are displayed in Table \ref{tab:tab2}.
\begin{table}[h!]
\centering
\begin{tabular}{ |C{1.3cm}||C{3.1cm}||C{2.2cm}|  }
 \hline
 \multicolumn{3}{|c|}{Design Parameters and Constraints} \\
 \hline
 Parameter & Baseline Value & Constraint \\
 \hline
$M$ & $64$ & $-$ \\
$N$ & $192$ & $-$ \\
$K$ & $24$ & $-$ \\
$L$ & $8$ & $-$ \\
$a_{mn}$ & Randomly distributed according to i.i.d. complex Normal distribution& Unbounded \\
%$k$ & $62.8319 \text{m}^{-1}$ & $k = 62.8319 \text{m}^{-1}$ \\
 \hline
\end{tabular}
\caption{Summary of design parameters and constraints for the general linear system sensing matrix design problem.}
\label{tab:tab2}
\end{table}

\begin{figure}[h!]
\centering
\includegraphics[width=8.375cm, clip=true]{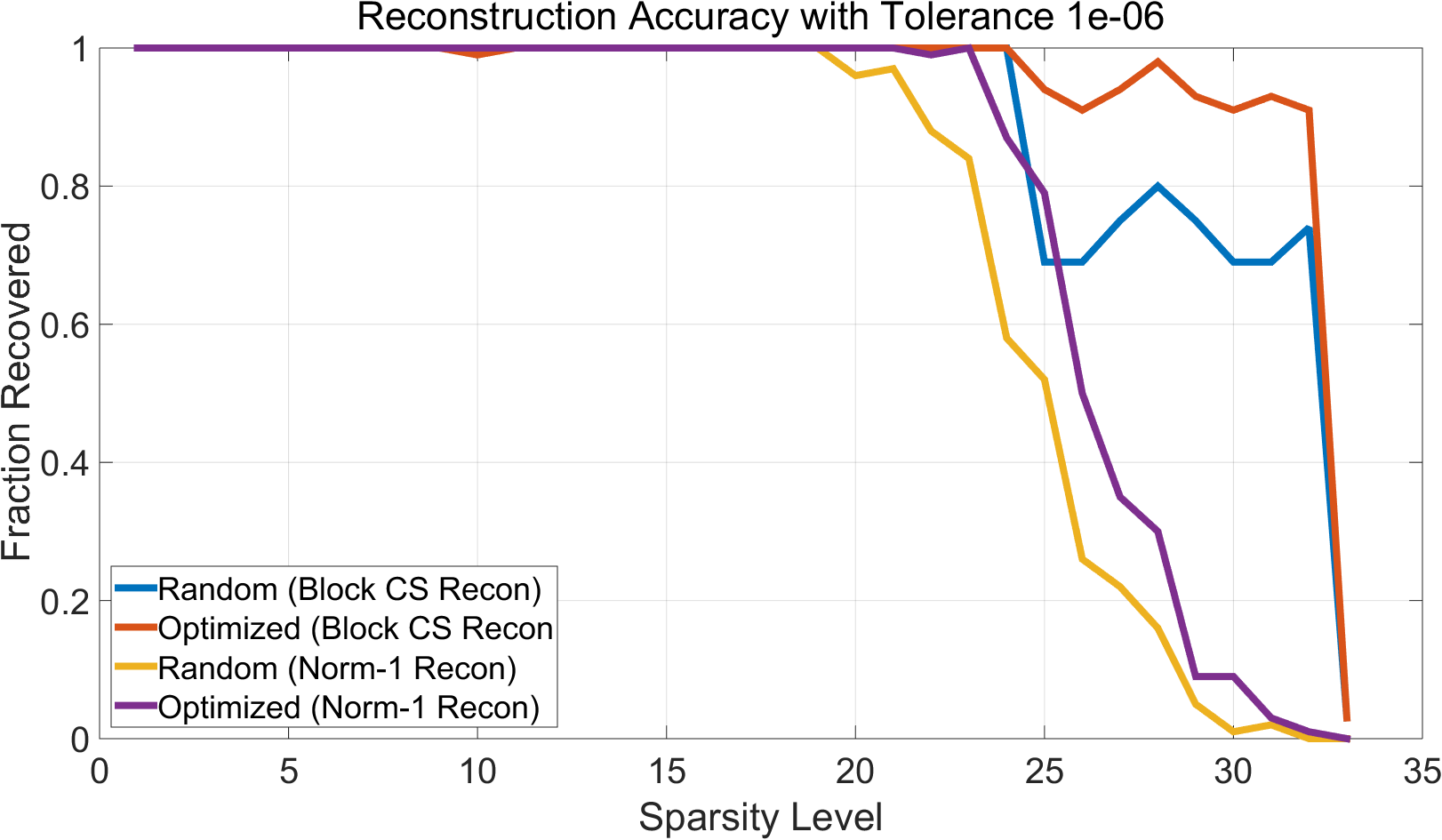}
    \caption{Numerical comparison of the reconstruction accuracies of joint $\ell_2/\ell_1$ reconstruction and standard $\ell_1$ reconstruction using the baseline random and optimized designs for the general linear system problem.}
   \label{fig:lin_results}
\end{figure}

For this design problem, the capacity of the optimized sensing matrix increased slightly, from $-2.7$ to $-0.8$. Nevertheless, the reconstruction accuracy was improved, as can be seen in Figure \ref{fig:lin_results}. Using joint $\ell_2/\ell_1$ minimization, the optimized sensing matrix reconstructed $> 90\%$ of block-sparse vectors up to a block sparsity $S_B = 4$ (total sparsity $S = 32$), whereas the baseline randomized sensing matrix only reconstructed $75-80\%$. Of the three examples presented in this paper, the reconstruction accuracy is unsurprisingly increased the least in this problem. Therefore, it may not be worth the effort to run the optimization procedure in applications, such as this one, where the designer has significant control over the sensing matrix.

%%%%%%%%%%%%%%%%%%%%%%%%%%%%%%%%%%%%%%%%%%%%%%%%%%%%%%%%%
\vspace{7pt}
\section{Conclusion}
\label{sec:conclusion}
This paper describes a novel method for designing sensing matrices with enhanced block-sparse signal recovery capabilities. By maximizing the minimum capacity over a set of sub-matrices selected from columns of the full sensing matrix, the design method is capable of significantly improving the reconstruction results obtained using joint $\ell_2/\ell_1$ minimization. This capability was demonstrated in three applications: a sparse pulse reconstruction problem, an electromagnetic imaging problem, and a general linear system. These results showed that the design method can be extremely beneficial in applications where the measurement system is constrained by practical limitations, but less beneficial when one has greater control over the sensing matrix.

\vspace{7pt}

\bibliographystyle{IEEEtran}
\bibliography{./SICA-TA}

\end{document}